\documentclass[oneside,11pt]{article}

\usepackage[utf8]{inputenc}
\usepackage[T1]{fontenc}

\usepackage{amsthm}
\usepackage{amsmath}
\usepackage{amsfonts}
\usepackage{amssymb}

\usepackage{graphicx} 
\usepackage{float}
\usepackage{booktabs}
\usepackage{multirow}
\usepackage{rotating}
\usepackage{array}
\usepackage{xcolor}
\usepackage[ruled]{algorithm2e}
\usepackage{subcaption}
\usepackage{ragged2e}
\usepackage{tikz}
\usetikzlibrary{angles, quotes, 3d, perspective}

\newcolumntype{L}[1]{>{\RaggedRight\arraybackslash}p{#1}}
\newcolumntype{R}[1]{>{\RaggedLeft\arraybackslash}p{#1}}
\newcolumntype{C}[1]{>{\centering\let\newline\\\arraybackslash\hspace{0pt}}m{#1}}

\usepackage{breakcites}

\usepackage{nowidow}

\usepackage{enumitem}

\usepackage[top=1.5cm,bottom=2cm,right=2.5cm,left=2.5cm]{geometry}
\usepackage{titling}

\usepackage{natbib}

\usepackage{ifplatform}

\usepackage{textcomp}

\ifwindows
\fi

\allowdisplaybreaks

\newcommand{\Sq}{\mathbb{S}^q}

\newcommand{\rd}{\mathrm{d}}

\newcommand{\bx}{\boldsymbol{x}}
\newcommand{\be}{\boldsymbol{e}}

\newcommand{\by}{\boldsymbol{y}}

\newcommand{\bX}{\boldsymbol{X}}

\newcommand{\bU}{\boldsymbol{U}}

\newcommand{\bR}{\boldsymbol{R}}
\newcommand{\bO}{\boldsymbol{O}}
\newcommand{\bmu}{\boldsymbol\mu}
\newcommand{\bu}{\boldsymbol{u}}

\newcommand{\bff}{\boldsymbol{f}}

\newcommand{\bbf}{\boldsymbol{f}}

\newcommand{\bGamma}{\boldsymbol\Gamma}

\newcommand{\Hcal}{\mathcal{H}}

\newcommand{\equald}{\stackrel{d}{=}}
\newcommand{\inlaw}{\rightsquigarrow}

\newcommand{\lrp}[1]{\left(#1\right)}

\newcommand{\om}[1]{\omega_{#1}}

\DeclareFontFamily{OT1}{pzc}{}
\DeclareFontShape{OT1}{pzc}{m}{it}{<-> s * [1.10] pzcmi7t}{}
\DeclareMathAlphabet{\mathpzc}{OT1}{pzc}{m}{it}

\newtheorem{theorem}{Theorem}

\newtheorem{proposition}{Proposition}

\usepackage[pdftex,final,bookmarksnumbered,bookmarksopen=false,breaklinks,colorlinks]{hyperref}
\hypersetup{
	final,
	bookmarksnumbered=true,
	bookmarksopen=true,
	bookmarksopenlevel=0,
	unicode=false,
	pdftoolbar=true,
	pdfmenubar=true,
	pdffitwindow=false,
	pdftitle={A stereographic test of spherical uniformity},
	pdfdisplaydoctitle=true,
	pdftoolbar=true,
	pdfmenubar=true,
	pdflang={English},
	pdfauthor={Alberto Fernandez-de-Marcos, Eduardo Garcia-Portugues},
	pdfsubject={arXiv paper},
	pdfcreator={Alberto Fernandez-de-Marcos},
	pdfproducer={Alberto Fernandez-de-Marcos},
	pdfkeywords={Spherical data}{Stereographic projection}{Uniformity tests}
	pdfnewwindow=true,
	breaklinks=true,
	hidelinks,       
	linkcolor=black,
	citecolor=black,
	filecolor=magenta,
	urlcolor=cyan
}

\newif\ifmain
\maintrue
\newif\ifsupplement
\supplementtrue

\newif\iffigstabs
\figstabstrue

\begin{document}

\ifmain

\title{A stereographic test of spherical uniformity}
\setlength{\droptitle}{-1cm}
\predate{}%
\postdate{}%
\date{}

\author{Alberto Fern\'andez-de-Marcos$^{1}$ and Eduardo Garc\'ia-Portugu\'es$^{1,2}$}
\footnotetext[1]{Department of Statistics, Universidad Carlos III de Madrid (Spain).}
\footnotetext[2]{Corresponding author. e-mail: \href{mailto:edgarcia@est-econ.uc3m.es}{edgarcia@est-econ.uc3m.es}.}
\maketitle

\begin{abstract}
	We introduce a test of uniformity for (hyper)spherical data motivated by the stereographic projection. The closed-form expression of the test statistic and its null asymptotic distribution are derived using Gegenbauer polynomials. The power against rotationally symmetric local alternatives is provided, and simulations illustrate the non-null asymptotic results. The stereographic test outperforms other tests in a testing scenario with antipodal dependence between observations.
\end{abstract}
\begin{flushleft}
\end{flushleft}

\section{Introduction}
\label{sec:intro}

One of the opening questions when dealing with data supported on the unit (hyper)sphere $\Sq:=\{\bx\in\mathbb{R}^{q+1}:\|\bx\|=1\}$, $q\geq1$, is whether it has a non-trivial distributional pattern, such as a preferred direction, or whether it is uniformly distributed on $\mathbb{S}^q$. Given a random sample $\bX_1,\ldots,\bX_n$ generated from an absolutely continuous distribution $\mathrm{P}$ on $\Sq$, this question is formalized as that of testing $\Hcal_0:\mathrm{P}=\mathrm{Unif}(\Sq)$ vs. $\Hcal_1:\mathrm{P}\neq\mathrm{Unif}(\Sq)$, where $\mathrm{Unif}(\Sq)$ is the uniform distribution on $\Sq$.

An extensive class of tests of uniformity on $\mathbb{S}^q$ is the class of \textit{Sobolev tests} first proposed in \cite{Beran1969}, \cite{Gine1975}, and \cite{Prentice1978}. The test statistics of this class are $V$-statistics of order two acting on the pairwise angles of the sample, $\{\theta_{ij}\}_{i,j=1}^n$, with $\theta_{ij}=\cos^{-1}(\bX_i'\bX_j)$. Some recently described properties of this class are given in \cite{Jammalamadaka2020}, where truncated and standardized Sobolev statistics are shown to have null normal limits, and \cite{Garcia-Portugues2023}, where the class of Sobolev tests is shown to be essentially bijective with that of the tests based on the empirical cumulative distribution function of the projected sample.

This article provides a new parameter-dependent family of tests of uniformity on $\Sq$, $q\geq2$. The new family is motivated by a pairwise stereographic projection in which each datum acts as the projection pole. This motivation is developed in Section \ref{sec:geom_motiv}, while Section \ref{sec:stat} presents the definition of the parametric family of test statistics. Asymptotic null distributions and powers against a class of local alternatives are studied in Sections \ref{sec:asymp} and \ref{sec:power}. Section \ref{sec:MC} illustrates the non-null asymptotic results through Monte Carlo experiments. Section \ref{sec:pow_example} gives a testing situation with antipodal dependence between observations where the presented family of tests is advantageous. Proofs can be found in the Appendix.

\section{Motivation}
\label{sec:geom_motiv}

Let $\bX_1,\ldots,\bX_n$ be a sample without ties on the sphere $\mathbb{S}^2$. The standard stereographic projection is the function $s_2:\mathbb{S}^2\backslash\{\be_3\}\to\mathbb{R}^2$ defined as $s_2(\bx) := (x_1/(1-x_3),x_2/(1-x_3))'$, with $\be_3=(0,0,1)'$ being the north pole. For each observation $\bX_i$, we define the \textit{stereographic projection of $\bX_i$ with respect to the pole $\bX_j$}, $j\neq i$, $s(\bX_i;\bX_j):=s_2(\bX^{\ast}_i)$, as the standard stereographic projection of $\bX^{\ast}_i:=\bR_j\bX_i$ after a rotation matrix $\bR_j$ such that $\bX^{\ast}_j:=\bR_j\bX_j=\be_3$ is applied; see Figure \ref{fig:stereo-proj} for geometric visualization. The triangle $\triangle \bO\, \bX_i^{\ast}\, \bX_j^{\ast}$, where $\bO$ is the origin, is isosceles. Consequently, the angle at vertex $s(\bX_i;\bX_j)$ in the right triangle $\triangle \bO\, s(\bX_i;\bX_j)\, \bX_j^{\ast}$ is $\theta_{ij}/2$. Then, by the definition of cotangent, $\| s(\bX_i;\bX_j) \| = \cot(\theta_{ij}/2)$. Although $\bR_{j}$ and, consequently, $s(\bX_i;\bX_j)$ are not unique, the norm $\| s(\bX_i;\bX_j) \|$ is well defined.

\begin{figure}[htpb!]
\vspace{-0.25cm}
\centering
\pgfmathsetmacro{\myaz}{10}
\pgfmathsetmacro{\angrad}{0.5}
\pgfmathsetmacro{\rotang}{87}
\begin{tikzpicture}[declare function={
    stereox(\x,\y)=2*\x/(1+\x*\x+\y*\y);
    stereoy(\x,\y)=2*\y/(1+\x*\x+\y*\y);
    stereoz(\x,\y)=(-1+\x*\x+\y*\y)/(1+\x*\x+\y*\y);
    angle(\x,\y)=atan(\y/\x);
    Px=1.7;Py=-1.3;
    amax=2.3;amin=2;},scale=2.25,
    line join=round,line cap=round,
    dot/.style={circle,fill,inner sep=1pt}]
\pgfdeclarelayer{background} 
\pgfdeclarelayer{foreground} 
\pgfsetlayers{background,main,foreground}
\path[use as bounding box] (-2,-1.1) rectangle (2,1.3); 
\path[save path=\pathSphere,ball color={rgb:red,1;green,1;blue,1}, fill opacity=0.25]
(0,0) circle[radius=1]; 
\begin{scope}[3d view={10}{15}] 
    \draw (-amin,amax) -- (-amin,-amax) coordinate (bl) -- (amax,-amax) 
    coordinate (br)-- (amax,amax)
    node[above left]{}; 
    \begin{scope}
        \draw (-amin,amax) -- (amax,amax); 
    \end{scope}
    \begin{scope}[canvas is xy plane at z=0]
        \draw[dashed] (\myaz:1) arc[start angle=\myaz,end angle=\myaz+180,radius=1];
        \draw (\myaz:1) arc[start angle=\myaz,end angle=\myaz-180,radius=1];
        \path[save path=\pathPlane] (\myaz:amax) -- (\myaz+180:amax) --(bl) -- (br) -- cycle;
        \begin{scope}
            \draw[use path=\pathSphere];
        \end{scope}
        \begin{pgfonlayer}{background}
            \fill[color={rgb:red,75;green,134;blue,198},fill opacity=0.2]
            (\myaz:1) arc[start angle=\myaz,end angle=\myaz-180,radius=1]
            -- (-amin,0) -- (-amin,amax) -- (amax,amax) -- (amax,0) -- cycle;
        \end{pgfonlayer}
        \fill[color={rgb:red,75;green,134;blue,198},fill opacity=0.2]
        (\myaz:1) arc[start angle=\myaz,end angle=\myaz-180,radius=1]
        -- (-amin,0) -- (-amin,-amax) -- (amax,-amax) -- (amax,0) -- cycle; 
    \end{scope}

    \draw (Px,Py,0) node[dot,label=below:{\small$s(\bX_i;\bX_j)$}](s){}
    -- ({stereox(Px,Py)},{stereoy(Px,Py)},{stereoz(Px,Py)}) 
    node[dot,label=right:{\small$\bX^{\ast}_i$}](Xj){};
    \draw[latex-] (0,0,1) node[dot,label=above:{\small$\bX^{\ast}_j$}](Xi){} -- (0,0,0)node[dot,label=left:{\small$\bO$}](O){}; 
    \draw[latex-] ({stereox(Px,Py)},{stereoy(Px,Py)},{stereoz(Px,Py)}) -- (0,0,0); 

    \draw[-latex, color=blue] (0,0,0) -- (0, {-sin(\rotang)}, {cos(\rotang)}) node[dot, color=blue, label=below:\small$\bX_j$](Xjr){};
    \draw[-latex, color=blue] (0,0,0) -- ({stereox(Px,Py)}, {stereoy(Px,Py)*cos(\rotang)-stereoz(Px,Py)*sin(\rotang)}, {stereoy(Px,Py)*sin(\rotang) + stereoz(Px,Py)*cos(\rotang)}) node[dot, color=blue, label=below:\small$\bX_i$](Xir){};
    
    \begin{scope}[rotate around z={angle(Px,Py)}, canvas is xz plane at y=0]
        \draw ({\angrad*cos(90-acos(stereoz(Px,Py)))},{\angrad*sin(90-acos(stereoz(Px,Py)))}) node[above, label={above:\small$\theta_{ij}$}]{} arc[start angle = {90 - acos(stereoz(Px,Py))}, end angle = 90, radius = \angrad];
        \draw ({(Px^2+Py^2)^0.5-\angrad*cos(acos(stereoz(Px,Py))/2)},{\angrad*sin(acos(stereoz(Px,Py))/2)}) node[below, label={\small$\frac{\theta_{ij}}{2}$}]{} arc[start angle = {180-acos(stereoz(Px,Py))/2}, end angle = 180, radius = \angrad];
    \end{scope}

    \begin{pgfonlayer}{background} 
        \draw[dashed] (Xj) -- (Xi);
    \end{pgfonlayer}

    \begin{scope}[rotate around x={\rotang}, rotate around z={angle(Px,Py)}, canvas is xz plane at y=0]
        \draw[dotted, color=blue] (0:1) arc[start angle=0,end angle=360,radius=1];
        \draw[color=blue] ({\angrad*cos(90-acos(stereoz(Px,Py)))},{\angrad*sin(90-acos(stereoz(Px,Py)))}) node[above, label={below:\small$\theta_{ij}$}]{} arc[start angle = {90 - acos(stereoz(Px,Py))}, end angle = 90, radius = \angrad];
    \end{scope}

    \begin{scope}[rotate around z={angle(Px,Py)}, canvas is xz plane at y=0]
        \draw[dotted] (0:1) arc[start angle=0,end angle=360,radius=1];
    \end{scope}

    \draw[line width=0.5mm, color = {rgb:red,230;green,40;blue,40}] (0,0,0) -- (s); 
    \draw (O) node[dot] {};
    \draw (s) node[dot] {};
    
\end{scope}
\end{tikzpicture}
\vspace{-0.25cm}
\caption{\small Stereographic projection $s(\bX_i;\bX_j)$ of $\bX_i$ with respect to the pole $\bX_j$. The points $\bX_i$ and $\bX_j$ have been rotated into $\bX_i^\ast$ and $\bX_j^\ast$ such that $\bX_j^\ast$ is the north pole. The dashed blue (black) circle represents the great circle $S^1$ that joins $\bX_i$ and $\bX_j$ ($\bX_i^\ast$ and $\bX_j^\ast$). The length of the red segment is the norm of the stereographic projection, $\|s(\bX_i;\bX_j)\|$, which is well defined.} \label{fig:stereo-proj}
\end{figure}
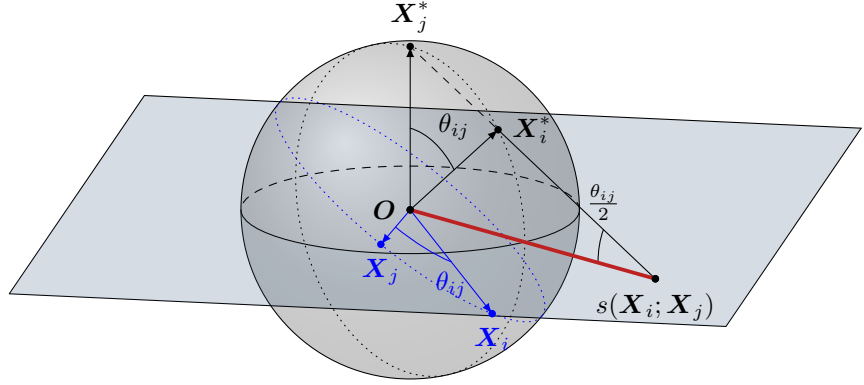

We consider the average of the lengths of the stereographic projections of pairwise observations as a test statistic:
\begin{align}
    P_n := \frac{2}{n}\sum_{1\leq i<j\leq n} \| s(\bX_i;\bX_j) \|. \label{eq:Sn}
\end{align}
This statistic is clearly invariant under rotations. A test based on it rejects the null hypothesis for large values of $P_n$, since $\|s(\bX_i;\bX_j)\|\to\infty$ as $\theta_{ij}\to 0$, which occurs under clustering of observations (indeed, the test statistic is not defined in the presence of ties in the sample). This observation is illustrated in Figure \ref{fig:stereo-distr}, which depicts the large summands of $P_n$ for samples drawn from a unimodal (central plot) and a bimodal distribution with antipodal symmetry (right plot). In contrast, if the sample is drawn from $\mathrm{Unif}(\mathbb{S}^2)$ (left plot), the summands of $P_n$ tend to be smaller than in the previous~cases.

\begin{figure}[!ht]
    \centering
    \begin{subfigure}[b]{0.33\textwidth}
        \centering
        \includegraphics[width=\textwidth,clip,trim={0cm 3cm 0cm 4cm}]{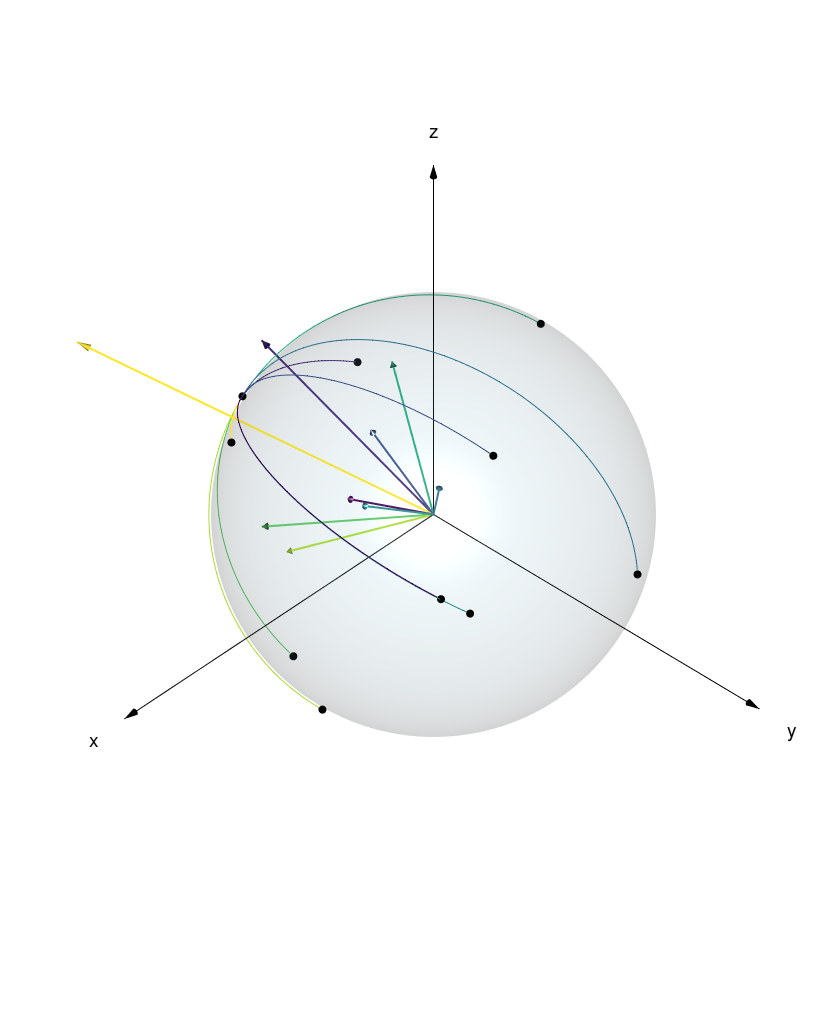}
        \vspace{-0.75cm}
    \end{subfigure}%
    \begin{subfigure}[b]{0.33\textwidth}
        \centering
        \includegraphics[width=\textwidth,clip,trim={0cm 3cm 0cm 4cm}]{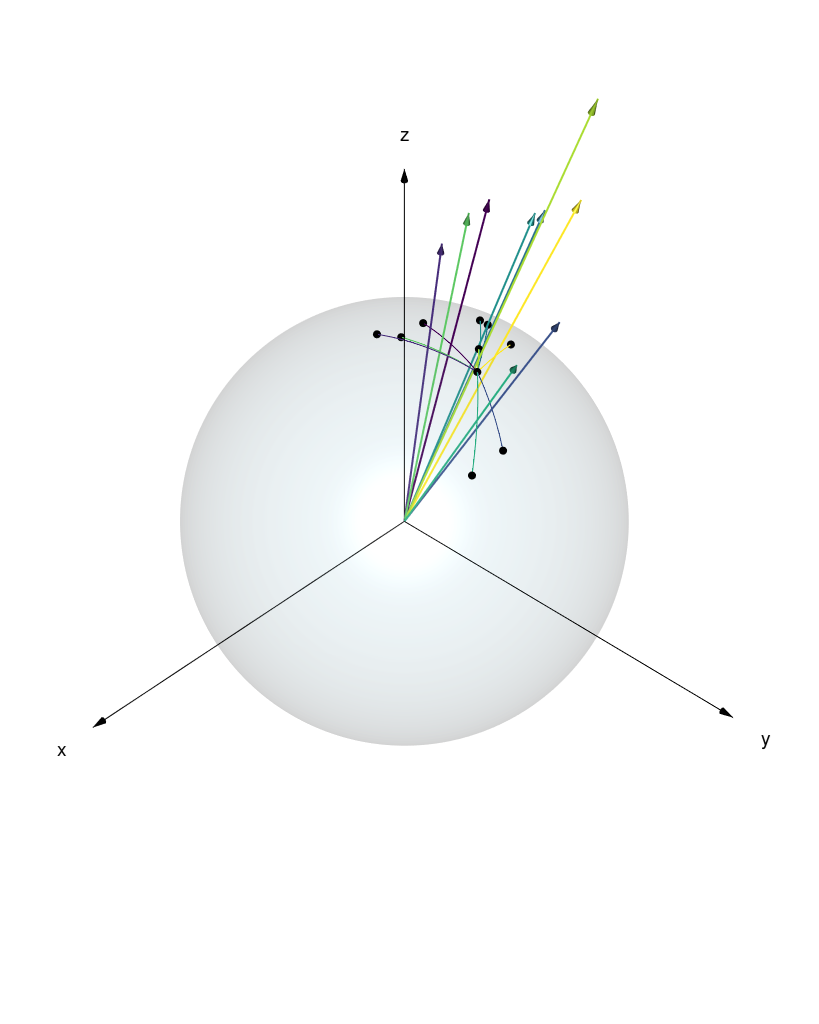}
        \vspace{-0.75cm}
    \end{subfigure}%
    \begin{subfigure}[b]{0.33\textwidth}
        \centering
        \includegraphics[width=\textwidth,clip,trim={0cm 3cm 0cm 4cm}]{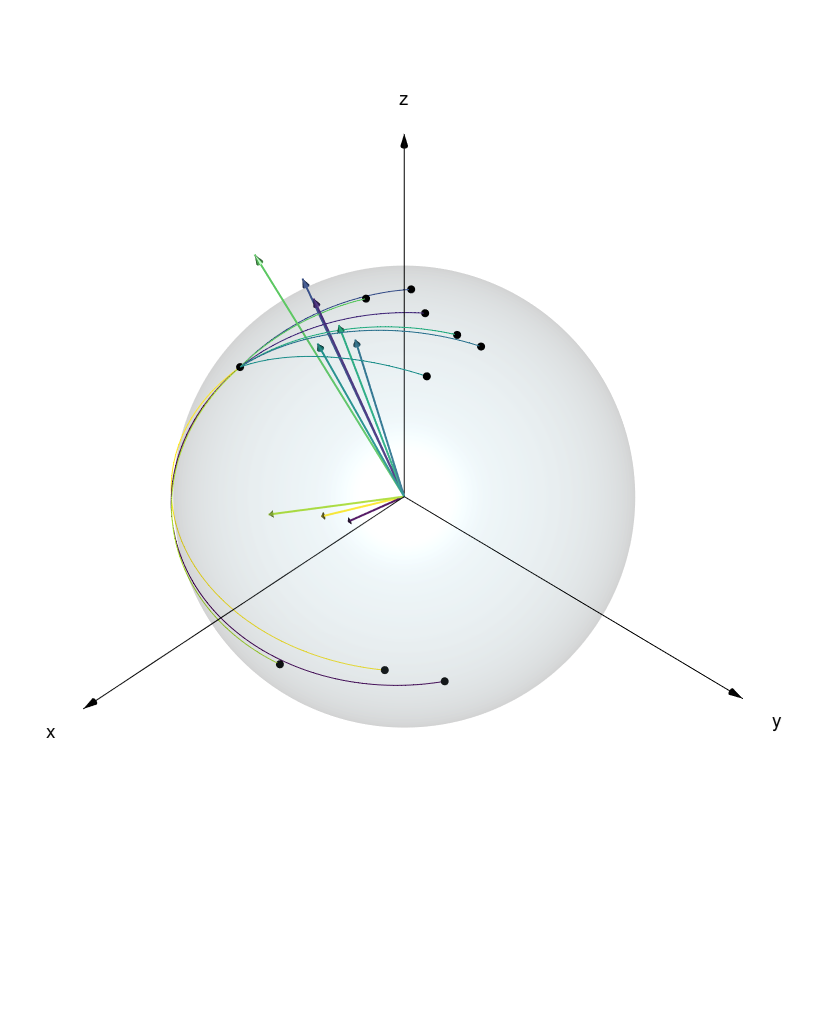}
        \vspace{-0.75cm}
    \end{subfigure}%
    \vspace*{-0.25cm}
    \caption{\small Illustration of the stereographic projections of $\bX_i$ with respect to the pole $\bX_1$, for $2\leq i\leq n=10$ on $\mathbb{S}^2$. The stereographic projection $s(\bX_i;\bX_1)$ is placed in the bisecting line of the angle $\theta_{i1}=\cos^{-1}(\mathbf{X}_i'\mathbf{X}_1)$, has a unique color that matches that of $\theta_{i1}$, and has the log-transformed norm $\log(1+\|s(\bX_i;\bX_1)\|)$ to improve visualization. The samples are generated from the uniform (left), von Mises--Fisher (center), and Watson (right) distributions. The latter distributions share location parameter $\bmu=\be_3$ and concentration $\kappa = 10$. Non-uniformity is detected with longer stereographic projections.
    }
    \label{fig:stereo-distr}
\end{figure}

The extension of $P_n$ to higher dimensions $q\geq2$ is straightforward. $P_n$ uses only the norm of $s(\bX_i;\bX_j)$, which is the same as that of the circular stereographic projection $\tilde{s}_1:S^1\backslash\{\be_3\}\subset\mathbb{S}^2\to\mathbb{R}$, with $S^1$ being the great circle that contains $\bX^{\ast}_i$ and $\bX^{\ast}_j=\be_3$; see Figure \ref{fig:stereo-proj}. $\tilde{s}_1$ is obtained from its standard version $\bx\in\mathbb{S}^1\mapsto s_1(\bx):= x_1/(1-x_2)$ after a parametrization of $S^1\cong \mathbb{S}^1$. Therefore, for any $q\geq 2$ and distinct $\bX_i$, $\bX_j$, we can define $s(\bX_i;\bX_j)$ such that $\|s(\bX_i;\bX_j)\|=|\tilde{s}_1(\bX_i^*)|$ using $\tilde{s}_1:S^1\backslash\{\be_{q+1}\}\subset\mathbb{S}^q\to\mathbb{R}$, i.e., the great circle $S^1\subset\mathbb{S}^q$ that joins $\bX^{\ast}_i$ and $\bX^{\ast}_j=\be_{q+1}$. As in the spherical case, the projection norm satisfies
\begin{align*}
    \|s(\bX_i;\bX_j)\|
    =\sqrt{\frac{1-((\bX_i^{\ast})_{q+1})^2}{(1-(\bX_i^{\ast})_{q+1})^2}}
    =\sqrt{\frac{1+\cos\theta_{ij}}{1-\cos\theta_{ij}}}
    =\cot(\theta_{ij}/2).
\end{align*}
This expression follows from $\smash{(\bX_i^{\ast})_{q+1}=\be_{q+1}^{\top}\bR_j\bX_i=(\bR_j\bX_j)^{\top}\bR_j\bX_i=\bX_j^{\top}\bX_i=\cos\theta_{ij}}$, $\|\bX_i^{\ast}\|=1$, and the stereographic projection formula. Since the projection norm depends only on $\bX_i^{\top}\bX_j$, it is independent of the choice of $\bR_j$, symmetric in $\bX_i$ and $\bX_j$, and invariant under a common rotation of both observations.

Then, for a sample $\bX_1,\ldots,\bX_n$ without ties on $\mathbb{S}^q$, $q\geq 2$, the stereographic projection statistic extending \eqref{eq:Sn} is defined as
\begin{align*}
    P_n := \frac{2}{n}\sum_{1\leq i<j\leq n} | \tilde{s}_1(\bR_j\bX_i) |=\frac{2}{n}\sum_{1\leq i<j\leq n} \cot(\theta_{ij}/2).
\end{align*}

\section{Test statistic}
\label{sec:stat}

Following the construction of the stereographic projection statistic $P_n$, we define a family of statistics for testing uniformity on $\mathbb{S}^q$, $q\geq2$, as
\begin{align}
    T_{n}(a):=\frac{2}{n}\sum_{1\leq i<j\leq n} \psi(\theta_{ij};a) - (n-1)\,\mathbb{E}_{\Hcal_0}\left[\psi(\theta_{12};a)\right],\label{eq:Tna}
\end{align}
where $\psi(\theta;a) := \cot(\theta/2) + a \tan(\theta/2)$, $a\in[-1,1]$, and
\begin{align}
    \mathbb{E}_{\Hcal_0}\left[\psi(\theta_{12};a)\right] = (1+a)\, \frac{(q-1)\Gamma((q-1)/2)^2}{2\Gamma(q/2)^2}\label{eq:EH0}
\end{align}
is obtained from Proposition \ref{prp:sob_coef} below. This family of statistics generalizes $P_n$, which is obtained for $a=0$: $T_n(0) = P_n - (n-1)\,\mathbb{E}_{\Hcal_0}\left[\psi(\theta_{12};0)\right]$. In addition, some statistics based on trigonometric kernels are particular cases of $T_n(a)$, such as $\psi(\theta; 1)=  2\csc(\theta)$ and $\psi(\theta; -1) = 2\cot(\theta)$. The election of $\theta\mapsto\tan(\theta/2)$ to extend $\psi(\theta;0)$ is justified because of (\textit{i}) the high similarity between the expressions of the Gegenbauer coefficients of both kernels and (\textit{ii}) their somewhat complementary behaviors. Indeed, $\tan(\theta/2)\to\infty$ as $\theta\to\pi$, and hence when $a\neq 0$, $T_n(a)$ also strongly penalizes antipodal observations, unlike $P_n$ and other classical Sobolev tests (see Section \ref{sec:pow_example}).

The statistic $T_n(a)$ is related to the general class of Sobolev statistics. When $q\geq2$, Sobolev test ($V$-)statistics can be defined as
\begin{align*}
    S_n(\{w_{k,q}\}):=\frac{1}{n}\sum_{i,j=1}^n \phi(\theta_{ij}),\quad \phi(\theta):=
    \sum_{k=1}^\infty \left(1+\frac{2k}{q-1}\right) w_{k,q}C_k^{(q-1)/2}(\cos\theta),
\end{align*}
where $\big\{C_k^{(q-1)/2}\big\}_{k=0}^\infty$ is the collection of Gegenbauer polynomials on $L_q^2([-1,1])$, the space of square-integrable real functions on $[-1,1]$ with respect to the weight $x\mapsto (1-x^2)^{q/2-1}$. Gegenbauer polynomials form an orthogonal basis on $L_q^2([-1,1])$: $\int_{-1}^1 C_k^{(q-1)/2}(x)C_\ell^{(q-1)/2}(x)(1-x^2)^{q/2-1}\,\mathrm{d}x=\delta_{k,\ell}c_{k,q}$ for $k,\ell\geq0$, where $\delta_{k,\ell}$ is the Kronecker delta, and
\begin{align}
    c_{k,q}:=\frac{\om{q}}{\om{q-1}}\lrp{1+\frac{2k}{q-1}}^{-2}d_{k,q},\quad d_{k,q}:=\lrp{1+\frac{2k}{q-1}}C_k^{(q-1)/2}(1), \label{eq:c_d_kq}
\end{align}
with $\omega_q := 2\pi^{(q+1)/2}/\Gamma((q + 1)/2)$ denoting the area of $\Sq$. The Gegenbauer coefficients of $\phi$, $\{b_{k, q}(\phi)\}_{k=0}^{\infty}$, are
\begin{align}
    b_{k, q}(\phi) := \frac{1}{c_{k,q}}\int_{-1}^{1} \phi(\cos^{-1}(x)) C_{k}^{(q-1)/2}(x)\,(1-x^2)^{q/2-1}  \, \mathrm{d} x\label{eq:sob_coef}
\end{align}
and correspond bijectively to the $\{w_{k, q}\}_{k=0}^{\infty}$ coefficients of its Sobolev form:
\begin{align}
    w_{k,q}=(1+2k/(q-1))^{-1}b_{k,q}(\phi).\label{eq:w_kq}
\end{align}

Taking $k=0$ in Proposition \ref{prp:sob_coef} gives \eqref{eq:EH0} from
\begin{align*}
    \mathbb{E}_{\Hcal_0}\left[\psi(\theta_{12};a)\right]&=\int_{\Sq} \psi(\cos^{-1}(\bx'\by);a)\,\rd\nu_q(\bx)\,\rd\nu_q(\by)\\
    &=\frac{\om{q-1}}{\om{q}}\int_{-1}^1 \psi(\cos^{-1}(x);a)(1-x^2)^{q/2-1}\,\rd x\\
    &=b_{0, q}(\psi(\cdot;a)),
\end{align*}
where $\nu_q$ represents the uniform measure on $\Sq$.

\begin{proposition}\label{prp:sob_coef}
Let $a\in[-1,1]$, $k\geq 0$, and $q\geq2$. Then,
\begin{align}
    b_{2k+1, q}(\psi(\cdot;a)) &= \alpha_{k,q}\frac{(2k+1)(4k+q+1)}{2k+q} \left(1 - a\right) ,\label{eq:sob_weight_odd}\\
    b_{2k, q}(\psi(\cdot;a)) &= \alpha_{k,q}(4k+q-1) \left(1 + a\right),\label{eq:sob_weight_even}
\end{align}
where $\alpha_{k,q} = \Gamma \left(k+1/2\right)^2\Gamma \left((q-1)/2\right)^2/[2 \pi\Gamma \left(k+q/2\right)^2]$.
\end{proposition}

The statistic $T_n(a)$ is not defined for the circular case $q=1$, where the stereographic projection is too extreme under uniformity: $\mathbb{E}_{\Hcal_0}\left[\psi(\theta_{12};a)\right]=\infty$ for any $a\in [-1,1]$. The sequence of Chebyshev coefficients $\{b_{k,1}(\psi(\cdot;a))\}_{k=0}^\infty$, $b_{k,1}(\psi(\cdot;a)):= c_{k,1}^{-1}\int_{-1}^{1} \psi(\cos^{-1}(x);a) T_{k}(x)\,(1-x^2)^{-1/2}\,\mathrm{d}x$, where $\{T_k\}_{k=0}^\infty$ are the Chebyshev polynomials, also has non-finite terms for any $a\in [-1,1]$.

For $q=2$, Section \ref{sec:asymp} shows that standard asymptotic results for degenerate $U$-statistics are not available for $T_{n}(a)$. Thus, we consider the \textit{$K(\geq1)$-truncated} version of the statistic $T_n(a)$ for $q\geq 2$:
\begin{align*}
    T_{n,K}(a) &:= \frac{2}{n}\sum_{1\leq i<j\leq n}\psi_{K}(\theta_{ij}; a) - (n-1)b_{0,q}(\psi(\cdot;a)),\\
    \psi_{K}(\theta; a)&:= \sum_{k=0}^{K} b_{k,q}(\psi(\cdot;a)) C_{k}^{(q-1)/2}(\cos \theta).
\end{align*}

\section{Null asymptotic distribution}
\label{sec:asymp}

Given an independent sample $\bX_1,\ldots,\bX_n$ from the distribution $\mathrm{P}$, the test based on $T_n(a)$ rejects $\Hcal_0:\mathrm{P}=\mathrm{Unif}(\mathbb{S}^q)$ at the asymptotic level $\alpha$ when $T_{n}(a)$ is larger than the upper-$\alpha$ quantile of its weak limit $T_{\infty}^{(0)}(a)$, given below.

\begin{theorem} \label{thm:asymp_distr}
    Let $a\in[-1,1]$ and $q\geq3$. Then, under $\Hcal_0$ and as $n\to\infty$, 
    \begin{align}
        T_{n}(a)\inlaw T_{\infty}^{(0)}(a)
        \equald\sum_{k=1}^\infty w_{k,q} (Y_k - d_{k,q}),\label{eq:null_asymp_dist}
    \end{align}
    with $Y_k\sim \chi^2_{d_{k,q}}$, $k\geq1$, being a sequence of independent random variables and $\equald$ denoting equality in distribution.
\end{theorem}

Since $w_{k,q}>0$ for $|a|<1$, the test based on $T_n(a)$ is consistent against any non-uniform density $f\in L^2(\Sq)$. Indeed, if $\bbf_k\in\mathbb{R}^{d_{k,q}}$ are the coefficients of $f$ on the $k$th-order subspace of spherical harmonics (orthogonal polynomials of degree $k$ in $L^2(\Sq)$), then $w_{k,q}>0$ allows detecting the non-uniformity of $f$ arising from $\|\bff_k\|>0$. When $a=1$ (respectively, $a=-1$), the test does not detect $\|\bff_k\|>0$ for all odd values (respectively, even) of~$k$.

Theorem \ref{thm:asymp_distr} does not apply for $q = 2$ since $\mathbb{E}_{\mathcal{H}_0}[\psi^2(\theta_{12};a)]=\|\psi(\cdot;a)\|^2_{L_q^2} \nless \infty$ and standard tools for degenerate $U$-statistics are not applicable. However, considering the $K$-truncated statistic $T_{n, K}(a)$, $\mathbb{E}_{\Hcal_0}[\psi^2_{K}(\cdot;a)]<\infty$ for $q\geq 2$. Then, as a particular case of Theorem \ref{thm:asymp_distr} where $w_{k,q}=0$ for $k>K$, $T_{n,K}(a)\inlaw T_{\infty,K}^{(0)}\equald\sum_{k=1}^{K}w_{k,q} (Y_k - d_{k,q})$ as $n\to\infty$ under $\Hcal_0$.

\section{Power against local alternatives}
\label{sec:power}

Let $\mathrm{P}_{\kappa,f}$ denote the distribution on $\Sq$ of the density function
\begin{align}
    \bx \mapsto \frac{c_{q, \kappa, f}}{\om{q-1}} f(\kappa\bx'\bmu),\label{eq:pdf_rotsym}
\end{align}
where $c_{q, \kappa, f}^{-1}:= \int_{-1}^{1} (1-s^2)^{q/2-1} f(s)\,\rd s$ is the normalizing constant. This density is rotationally symmetric about $\bmu\in\mathbb{S}^q$, which is set without loss of generality to $\bmu = \be_{q+1}$. The non-uniformity of \eqref{eq:pdf_rotsym} is controlled by the \emph{concentration} $\kappa>0$ and the \emph{angular function} $f:\mathbb{R}\rightarrow\mathbb{R}_0^{+}$, where $f(0) = 1$. If $\kappa=0$ or $f\equiv1$, then $\mathrm{P}_{\kappa,f}=\mathrm{Unif}(\mathbb{S}^d)$. We consider a sequence of local alternatives to uniformity letting $\kappa_n\to 0$ as $n\to \infty$ at a suitable rate, for a fixed $f$. For this purpose, the distribution of the independent and identically distributed (iid) sample $\bX_1,\ldots,\bX_n$ generated from \eqref{eq:pdf_rotsym} with varying $\kappa_n$ is denoted $\mathrm{P}_{\kappa_n,f}^{(n)}$.

The following result, which follows from Theorem 3 in \cite{Garcia-Portugues:Sobolev}, gives the non-null behavior of $T_{n,K}(a)$ under $\mathrm{P}_{\kappa_n,f}^{(n)}$.

\begin{theorem} \label{thm:power}
Let $a\in[-1,1]$, $q \geq 2$, and $f$ be $K$ times differentiable at zero. Define $k_v:=1 + \mathbb{I}[a=1]$ and fix $\tau\in\mathbb{R}$. Then, assuming that $K\geq k_v$, under $\mathrm{P}_{\kappa_n,f}^{(n)}$ with $\kappa_n=n^{-1/(2k_v)}\tau$,
\begin{align}
T_{n,K}(a)\inlaw T_{\infty,K}^{(1)}(a)\equald w_{k_v,q}(\widetilde{Y}_{k_v}-d_{k_v,q})+\sum_{k=1,\,k\neq k_v}^{K}w_{k,q}\,(Y_k-d_{k,q}),\label{eq:nonnull_asymp_dist}
\end{align}
where $\widetilde{Y}_{k}\sim \chi^2_{d_{k,q}}(\xi_{k,q}(\tau))$ and $\xi_{k,q}(\tau):=d_{k,q}(f^{\underline{k}}(0))^2\tau^{2k}{\prod_{\ell=0}^{k-1}(2\ell+q+1)^{-2}}$ ($f^{\underline{k}}$ denotes the $k$th derivative of $f$).
\end{theorem}

Under $\mathrm{P}_{\kappa_n,f}^{(n)}$ and provided that $f^{\underline{k_v}}(0)\neq0$, $\kappa_n\sim n^{-1/(2k_v)}$ is the \emph{detection threshold} of the test based on $T_{n,K}(a)$. That is, when $n^{1/(2k_v)}\kappa_n\to0$, $T_{n,K}(a)\inlaw T_{\infty,K}^{(0)}$ and, therefore, the test is blind to these alternatives; and when $n^{1/(2k_v)}\kappa_n\to\infty$, $\mathrm{P}_{\kappa_n,f}^{(n)}\left[T_{n,K}(a)>c\right]\to1$ for any $c>0$, so the test is consistent. The detection thresholds when $f^{\underline{k_v}}(0)=0$ are more convoluted; see Theorem 3 in \cite{Garcia-Portugues:Sobolev}.

The asymptotic power against local alternatives and the detection thresholds of the tests based on $T_{n}(a)$ are conjectured to be analogous to those in Theorem \ref{thm:power} for $T_{n,K}(a)$, letting $K\to\infty$. However, proving this claim is technical and involves adapting the proof of Theorem 4 in \cite{Garcia-Portugues:Sobolev} to Sobolev $U$-statistics, which is beyond the scope of this work. Monte Carlo experiments will be used to empirically evaluate this claim.

\section{Local power experiments}
\label{sec:MC}

To illustrate the theoretical behavior of the tests based on $T_{n,K}(a)$ and $T_n(a)$ under local alternatives $\mathrm{P}_{\kappa_n,f}^{(n)}$, and following \cite{Garcia-Portugues:Sobolev}, Monte Carlo experiments were carried out under three types of alternatives: the von Mises--Fisher (vMF) distribution with angular function $f_1(s)=\exp(s)$, the axial-symmetric distribution with $f_2(s)=\cosh(s)$ (an equal mixture of two vMFs), and the small circle distribution with $f_3(s;\nu)=\exp(-(s-\nu)^2)$, where $\nu=0.25$ determines the modal strip of the distribution. $M=10^4$ samples of size $n\in\{500, 2500\}$ were drawn from each of the three distributions on $\mathbb{S}^2$ (Figure \ref{fig:power-local-p3}) and $\mathbb{S}^3$ (Figure \ref{fig:power-local-p4}) with concentration $\kappa_n=n^{-1/\ell}\tau$ for each $\ell\in\{2,4,6\}$ and each $\tau=\{0,0.5,1,\ldots,6\}$. The tests based on $T_{n,K}(a)$ (for $\mathbb{S}^2$) and $T_n(a)$ (for $\mathbb{S}^3$) with $a\in\{-1,0,1\}$, and the Rayleigh and Bingham tests were performed at the asymptotic level $\alpha=5\%$. Simulations were performed using the \texttt{sphunif} R package \citep{Garcia-Portugues:sphunif}, where the new \texttt{Stereo} test is available.

\begin{figure}[!t]
    \centering
    \begin{subfigure}{\textwidth}
        \vspace*{-0.45cm}
        \includegraphics[width=\textwidth]{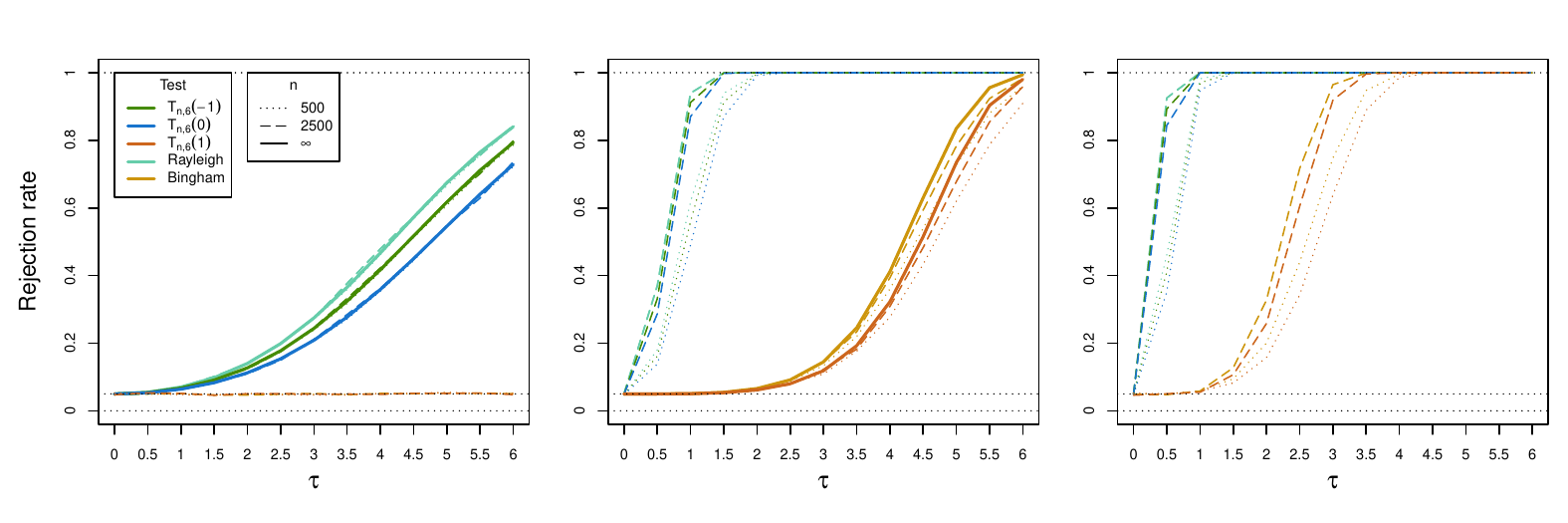}\vspace*{-0.4cm}
        \caption{\small $f_1(s)=\exp(s)$}
    \end{subfigure}
    \begin{subfigure}{\textwidth}
        \vspace*{-0.45cm}
        \includegraphics[width=\textwidth]{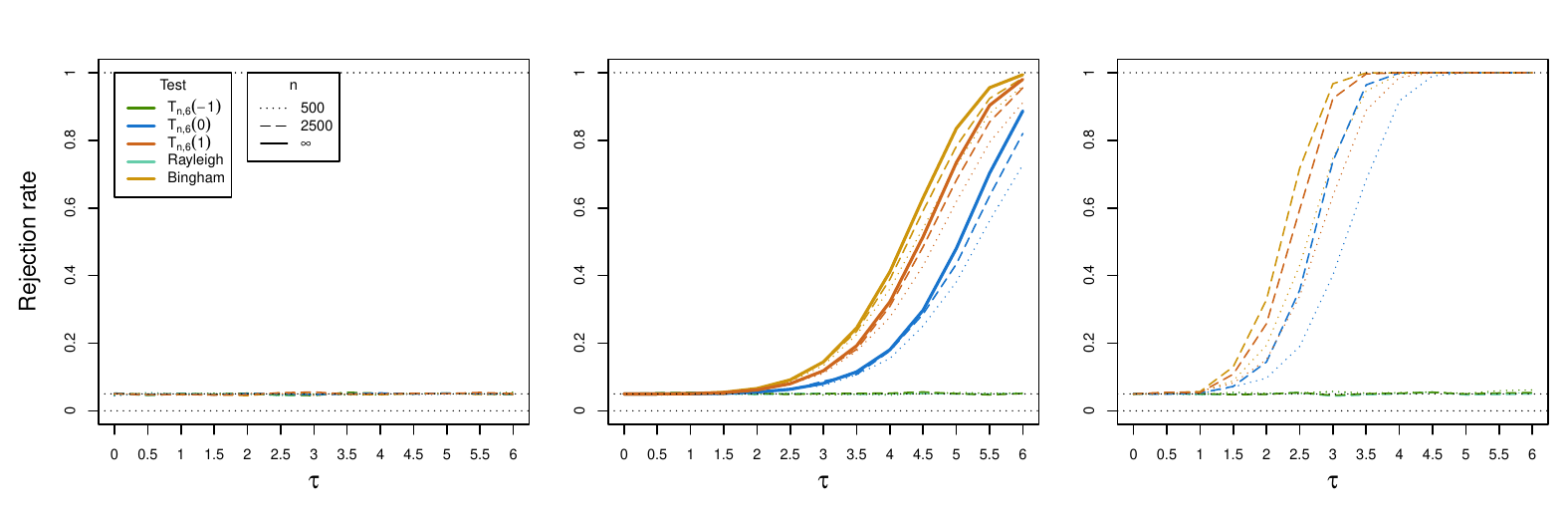}\vspace*{-0.4cm}
        \caption{\small $f_2(s)=\cosh(s)$}
    \end{subfigure}
    \begin{subfigure}{\textwidth}
        \vspace*{-0.45cm}
        \includegraphics[width=\textwidth]{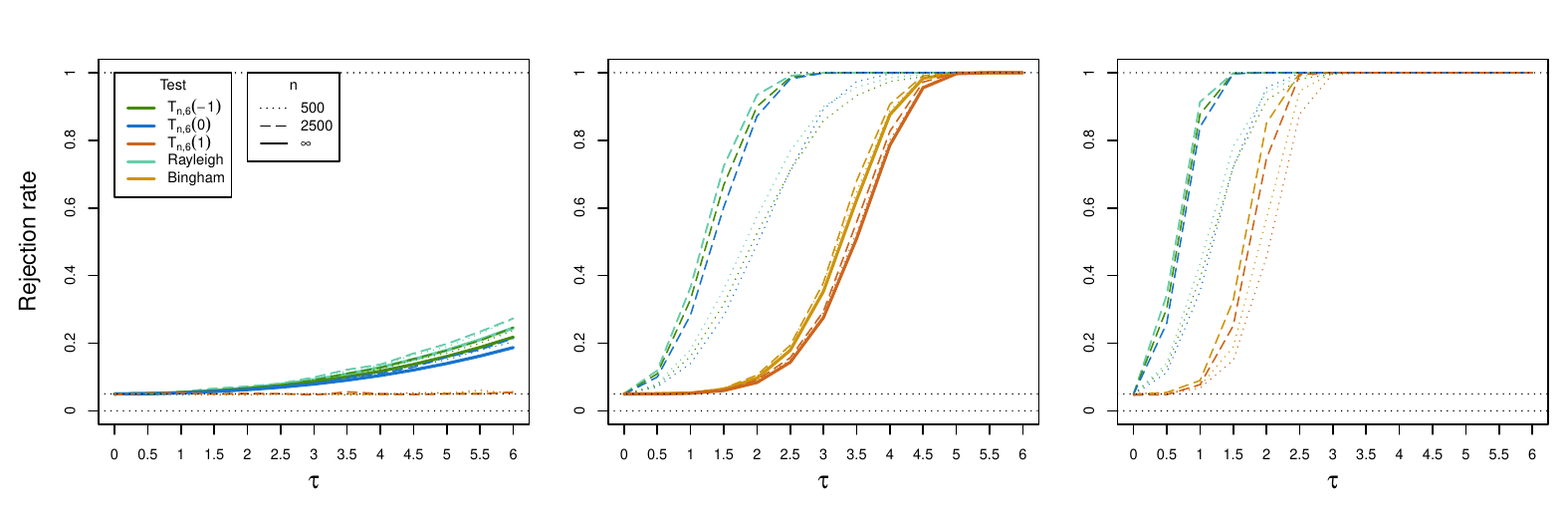}\vspace*{-0.4cm}
        \caption{\small $f_3(s;\nu)=\exp(-(s-\nu)^2)$}
    \end{subfigure}
    \caption{\small Rejection rates of the truncated test based on $T_{n,6}(a)$ for $a\in\{-1,0,1\}$, the Rayleigh test, and the Bingham test. Samples of sizes $n=500,2500$ (dotted and dashed lines) were drawn from local alternatives (a)--(c) on $\mathbb{S}^2$ with concentrations $\kappa_n=n^{-1/\ell}\tau$. From left to right, columns represent $\ell = 2,4,6$. Asymptotic powers are provided in solid lines when they are strictly between $\alpha$ and $1$. The asymptotic level is $\alpha=5\%$.}
    \label{fig:power-local-p3}
\end{figure}

\begin{figure}[!ht]
    \centering
    \begin{subfigure}{\textwidth}
        \vspace*{-0.45cm}
        \includegraphics[width=\textwidth]{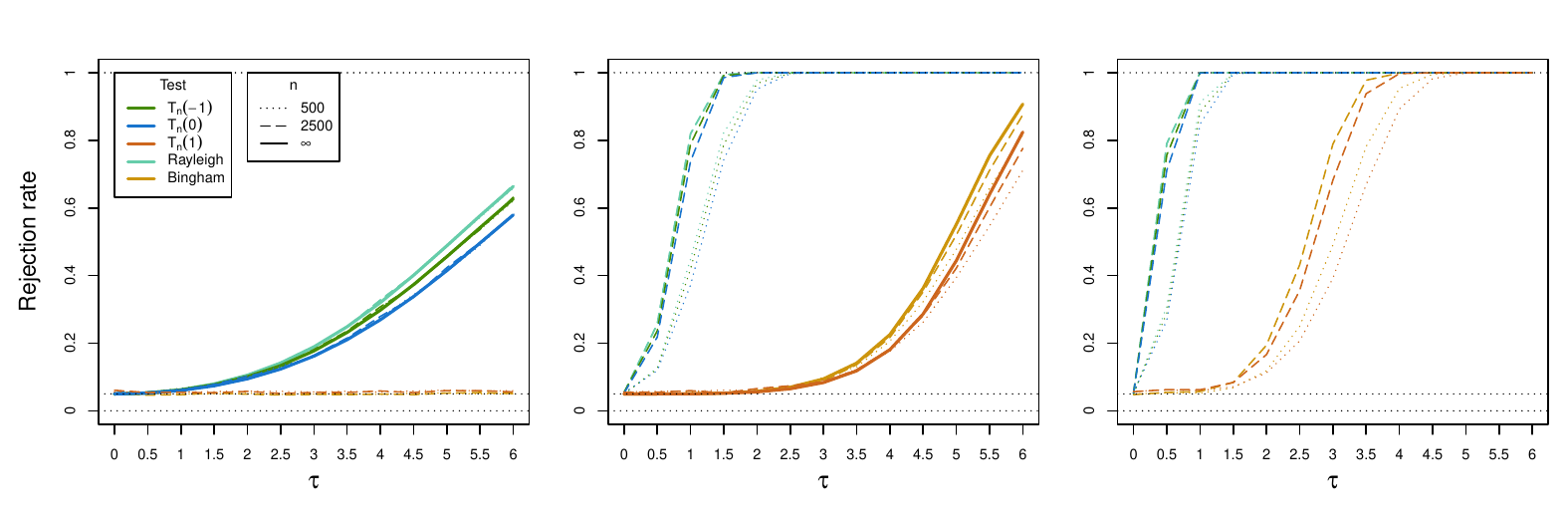}\vspace*{-0.4cm}
        \caption{\small $f_1(s)=\exp(s)$}\label{fig:power-local-vMF}
    \end{subfigure}
    \begin{subfigure}{\textwidth}
        \vspace*{-0.45cm}
        \includegraphics[width=\textwidth]{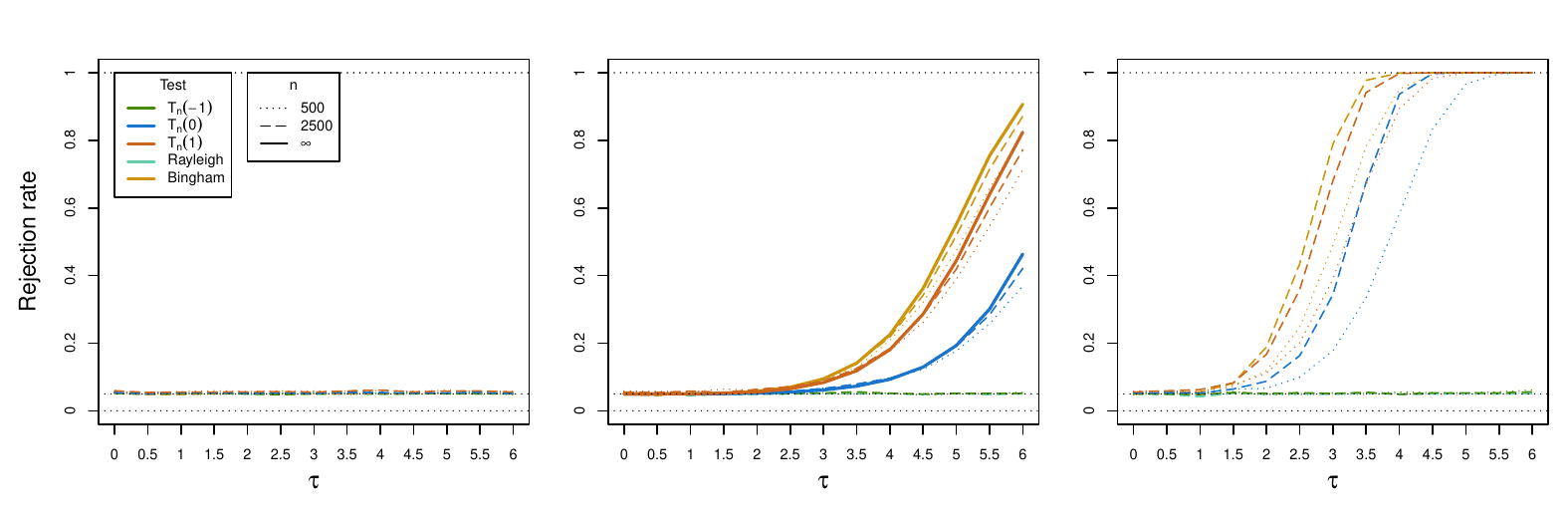}\vspace*{-0.4cm}
        \caption{\small $f_2(s)=\cosh(s)$}\label{fig:power-local-W}
    \end{subfigure}
    \begin{subfigure}{\textwidth}
        \vspace*{-0.45cm}
        \includegraphics[width=\textwidth]{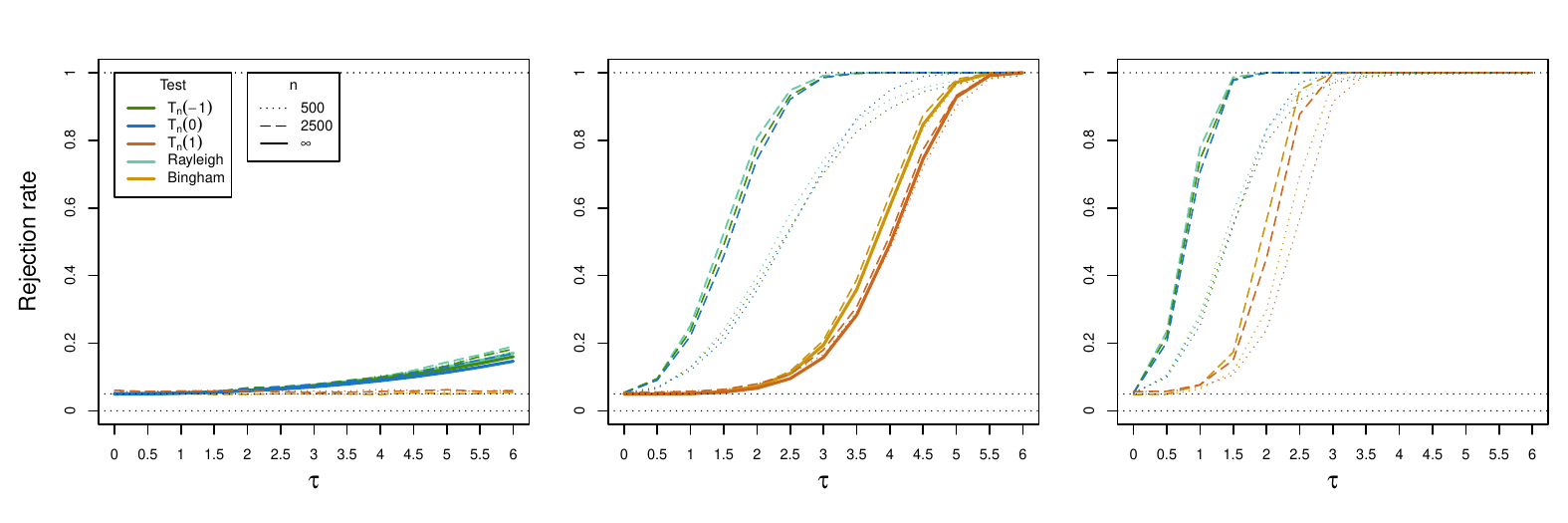}\vspace*{-0.4cm}
        \caption{\small $f_3(s;\nu)=\exp(-(s-\nu)^2)$}\label{fig:power-local-SC}
    \end{subfigure}
    \caption{\small Rejection rates of the test based on $T_n(a)$ for $a\in\{-1,0,1\}$, the Rayleigh test, and the Bingham test. The description of the figure is analogous to Figure \ref{fig:power-local-p3}, now with the local alternatives defined on $\mathbb{S}^3$.}
    \label{fig:power-local-p4}
\end{figure}

In Figures \ref{fig:power-local-p3} and \ref{fig:power-local-p4} the empirical rejection rates of the tests are plotted as a function of $\tau$. Asymptotic powers are provided when they are strictly between the level $\alpha$ and $1$, that is, when $\kappa_n$ coincides with each test's detection threshold depending on the local alternative. The detection thresholds of the tests based on $T_{n,K}(a)$ (from Theorem \ref{thm:power}) are as follows:
\begin{itemize}
    \item $\kappa_n\sim n^{-1/4}$ (center column) when $a=1$ (orange curves) for every alternative considered. This is the same detection threshold as the Bingham test (yellow).
    \item $\kappa_n\sim n^{-1/2}$ (left column) when $a<1$ (blue and green) for the vMF (top row) and the small circle (bottom row) alternatives. This is the same detection threshold as the Rayleigh test (cyan).
    \item $\kappa_n\sim n^{-1/4}$ (center column) when $a>-1$ (blue and orange) for the vMF mixture alternative (center row). The stereographic test when $a=-1$ (green) and the Rayleigh test (cyan) are blind to this alternative. This threshold is conjectured from Theorem 3 in \cite{Garcia-Portugues:Sobolev} (case $f^{\underline{k_v}}(0)=0$).
\end{itemize}
The results for $\mathbb{S}^3$ suggest that the detection thresholds also hold for $T_n(a)$.
 
\section{An advantageous testing situation}
\label{sec:pow_example}

The flexibility of $\psi(\theta;a)$ allows the proposed family of statistics to show a variety of behaviors. If $a=0$, the stereographic projection statistic $P_n$ arises, and nearby pairs of observations are strongly penalized since $\cot(\theta/2)\to\infty$ as $\theta\to 0$. When $a=1$, $\psi(\theta;1)=\cot(\theta/2)+\tan(\theta/2)$ results in a combination of the behavior shown by $P_n$ for close-by observations and a heavy penalization for pairs of observations that are almost antipodal, since $\tan(\theta/2)\to\infty$ as $\theta\to \pi$. This characteristic motivates a \textit{Uniform Antipodal-Dependent} (UAD) data-generating process where the antipodal dependence between observations is only well detected by $T_n(a)$ with $a\approx1$.

We consider the uniform distribution on a spherical cap centered at $\bmu\in\Sq$ with an opening angle $\theta\in[0,\pi]$ about $\bmu$, and denote it by $\mathrm{Unif}(\Sq; \bmu, \theta)$. The density of $\mathrm{Unif}(\Sq; \bmu, \theta)$ is $\bx \mapsto c_{q,\theta}\,1_{[1-r, 1]}(\bx'\bmu)$, where $r:=\cos(\theta)$ and $c_{q,\theta}$ is a normalizing constant. Since $\bX\sim\mathrm{Unif}(\Sq; \bmu, \theta)$ is a rotationally symmetric distribution about $\bmu$, the tangent-normal decomposition gives $\bX \equald V_{\bmu}(\bX) \bmu + (1-V_{\bmu}(\bX)^2)^{1/2} \bGamma_{\bmu}\bU_{\bmu}(\bX)$, where $\bGamma_{\bmu}$ is a $p\times(p-1)$ matrix whose columns are an orthogonal complement to $\bmu$, $v_{\bmu}(\bx):=\bx'\bmu$, $\bu_{\bmu}(\bx):=\bGamma_{\bmu}\bx/(1-v_{\bmu}(\bx)^2)^{1/2}$, and $\bu_{\bmu}(\bX)\sim \mathrm{Unif}(\mathbb{S}^{q-1})$ is independent from $v_{\bmu}(\bX)$ (see, e.g., Section 2 in \cite{Garcia-Portugues2020}). The density of $V_{\bmu}(\bX)$ is $\tilde{g}_q(v):=\omega_{q-1}c_{q,\theta}(1-v^2)^{q/2-1}1_{[1-r,1]}(v),$
where $c_{q,\theta}^{-1} :=\om{q-1}\int_{1-r}^{1} (1-t^2)^{q/2-1} \,\rd t = \om{q}[1-F_q(1-r)],$ with $F_q$ denoting the projected uniform distribution (see, e.g., Section 2.2 in \cite{Garcia-Portugues2023}). 
Then, the cumulative distribution of $V_{\bmu}(\bX)$ is $\tilde{G}_q(v)=[1-F_q(1-r)]^{-1} \{F_q(v)-F_q(1-r)\}$ and the quantile function is $\tilde{G}_q^{-1}(u)=F_q^{-1}\lrp{[1-F_q(1-r)]u+F_q(1-r)}$. Exact sampling for $\mathrm{Unif}(\Sq; \bmu, \theta)$ follows by the inversion method and the tangent-normal representation.

The UAD data-generating process of a sample $\bX_1,\ldots,\bX_n$ on $\Sq$ consists in $\lceil n/2\rceil$ iid observations drawn uniformly from the hypersphere, i.e. $\bX_i\sim\mathrm{Unif}(\Sq)$, for $1\leq i\leq \lceil n/2\rceil$, and the remaining $\lfloor n/2\rfloor$ observations generated following $\bX_{i+\lceil n/2\rceil}|\bX_i\sim \mathrm{Unif}(\Sq;-\bX_i, \theta)$, for $1 \leq i \leq \lfloor n/2\rfloor$. This process generates an identically distributed, but not independent, sample. Figure \ref{fig:antipodal-dist} shows three samples on $\mathbb{S}^2$ generated through this process for increasing spherical cap angles, suggesting how difficult it is to visually determine whether the sample comes from a uniform distribution or not.

\begin{figure}[!ht]
    \centering
    \begin{subfigure}[b]{0.25\textwidth}
        \centering
        \includegraphics[width=\textwidth,clip,trim={0cm 0cm 0cm 0cm}]{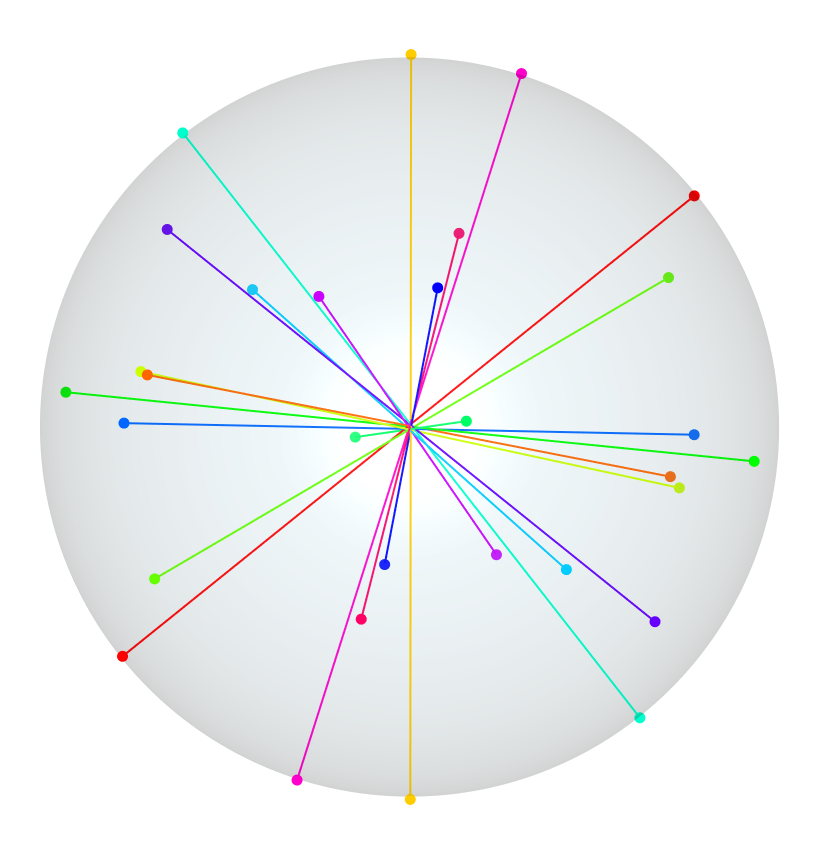}
        \vspace{-0.75cm}
    \end{subfigure}%
    \begin{subfigure}[b]{0.25\textwidth}
        \centering
        \includegraphics[width=\textwidth,clip,trim={0cm 0cm 0cm 0cm}]{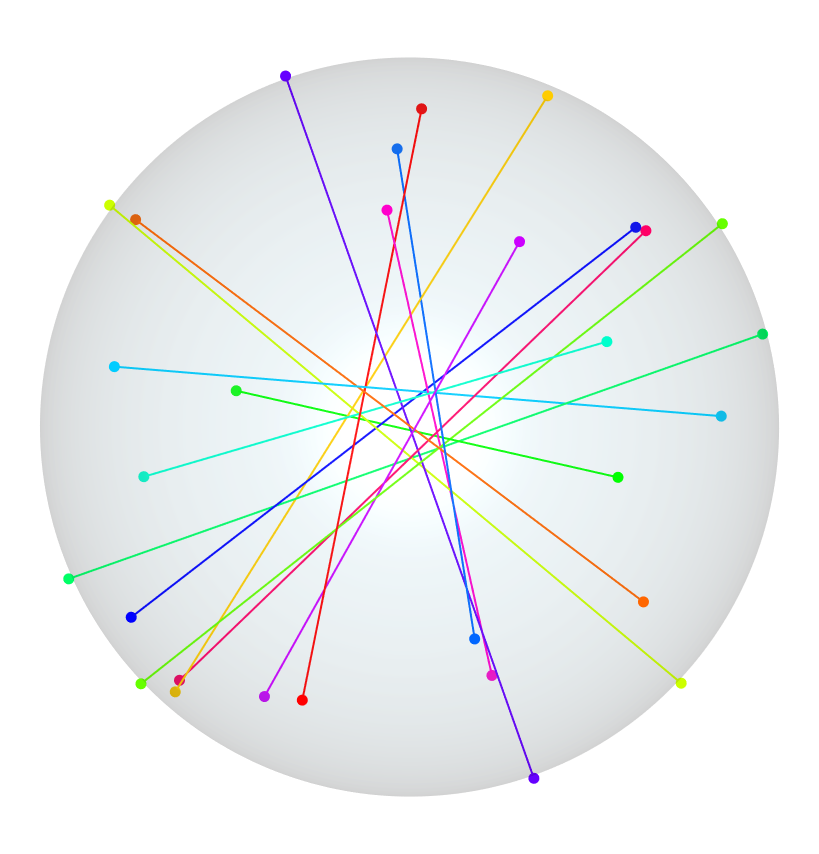}
        \vspace{-0.75cm}
    \end{subfigure}%
    \begin{subfigure}[b]{0.25\textwidth}
        \centering
        \includegraphics[width=\textwidth,clip,trim={0cm 0cm 0cm 0cm}]{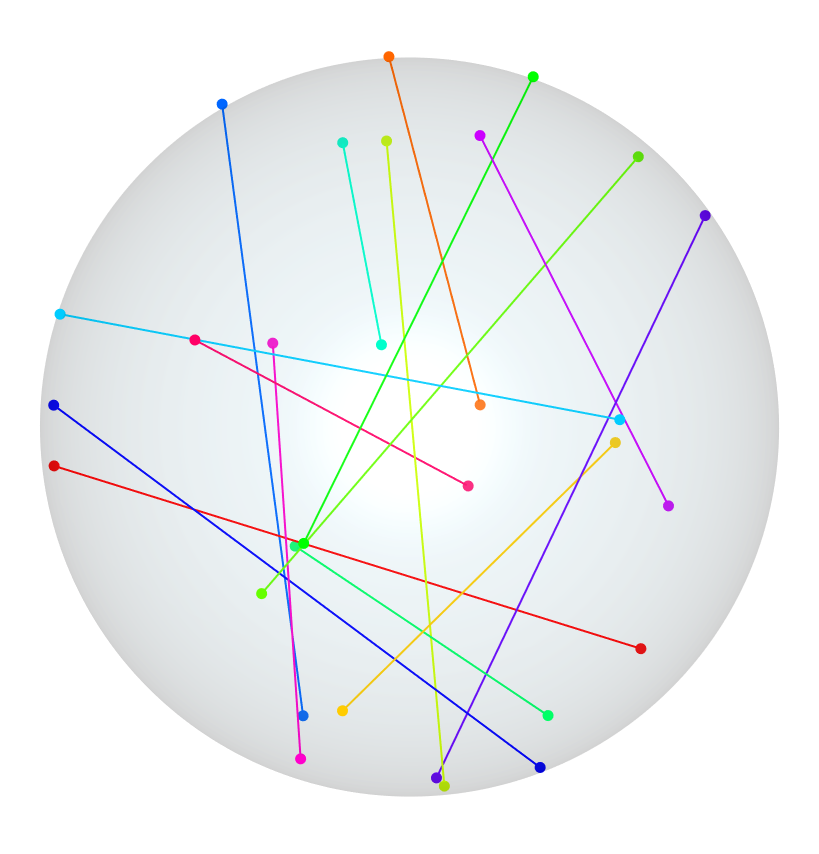}
        \vspace{-0.75cm}
    \end{subfigure}%
    \begin{subfigure}[b]{0.25\textwidth}
        \centering
        \includegraphics[width=\textwidth,clip,trim={0cm 0cm 0cm 0cm}]{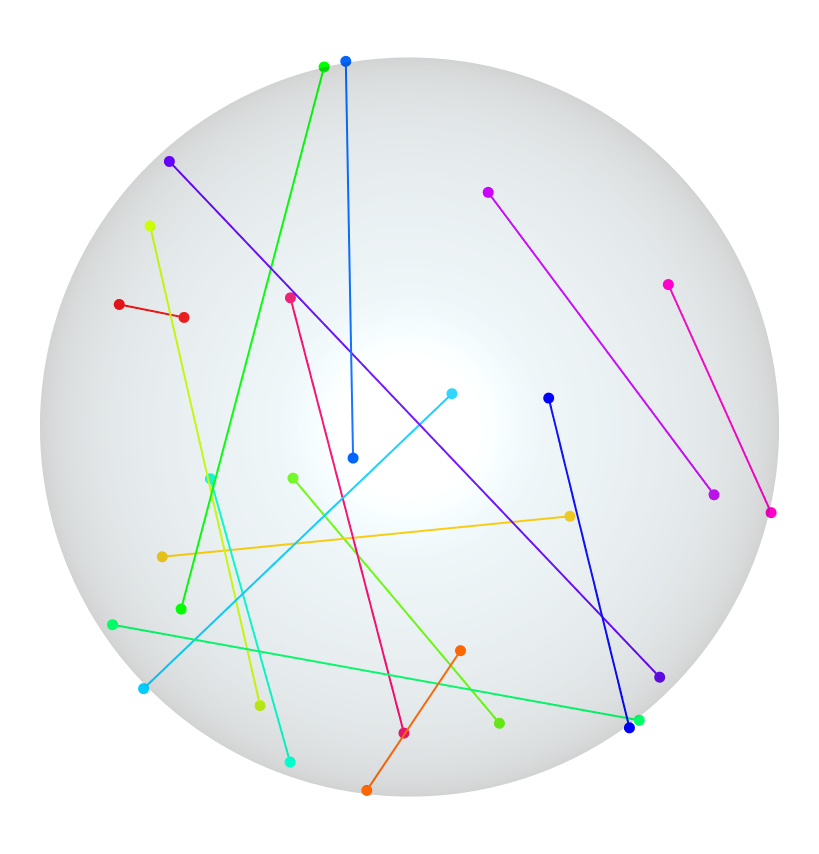}
        \vspace{-0.75cm}
    \end{subfigure}
    \caption{\small A sample of size $n=30$ generated by the UAD process on $\mathbb{S}^2$ with angle $\theta$ for $\theta\in\{1^\circ, 20^\circ, 90^\circ, 180^{\circ}\}$, from left to right, respectively. Each dependent pair is linked through a segment of the same color. The case $\theta=180^\circ$ is an iid sample from $\mathrm{Unif}(\mathbb{S}^2)$.}
    \label{fig:antipodal-dist}
\end{figure}

The automatic selection of $a$ for a given sample is possible using the $K$-fold test methodology from Section 5 in \cite{Fernandez-de-Marcos2023b}, which selects $a$ by a cross-validation procedure that maximizes the power under the learned testing scenario. A minor adaptation is required for $q=2$, as the variance of $T_n(a)$ under $\Hcal_0$ is non-finite for every parameter $a$: the optimal parameter $\hat{a}(\mathcal{S}_k)$ selected in Step (\textit{ii-a}) of Definition 2 in \cite{Fernandez-de-Marcos2023b} is set to $\hat{a}(\mathcal{S}_k)=\arg\max_{a\in\mathcal{A}} T_{|\mathcal{S}_k|}(a)$. We denote the resulting $K$-fold test statistic by $T_n^{(K)}$.

We compare the empirical performance of the tests based on $T_n(a)$ for $a\in\{-1,0,1\}$, the $10$-fold statistic $T_n^{(10)}$, the $10$-fold Poisson$^{(10)}$ and Softmax$^{(10)}$ statistics \citep{Fernandez-de-Marcos2023b}, and the well-known Rayleigh and Bingham statistics. For dimensions $q=2,3$ and spherical cap angles $\theta=1^\circ,10^\circ,20^\circ,45^\circ,90^\circ,135^\circ,180^\circ$, we generated $M=10^4$ Monte Carlo samples of size $n=100$ from the UAD process. The compared tests reject $\Hcal_0$ based on exact-$n$ critical values approximated by $M = 10^5$ Monte Carlo replicates at the significance level $\alpha = 5\%$. Table~\ref{tbl:unif-antipodal} shows the empirical rejection frequencies. When there is a high degree of antipodal dependence ($\theta< 20^\circ$ in $\mathbb{S}^2$ and $\theta\leq 20^\circ$ in $\mathbb{S}^3$), $T_n(1)$ shows a higher power than the rest of the tests. However, for $a=-1,0$, the tests based on $T_n(a)$ are blind, showing rejection rates similar to Poisson$^{(10)}$ and Softmax$^{(10)}$ for any value of $\theta$. Rayleigh is completely blind for every $\theta$, having rejection rates below $\alpha$ and even reaching $0\%$, due to the antipodal dependence. Bingham maintains a moderate rejection rate when $\theta<45^\circ$, since it is designed to capture axial alternatives, achieving a rejection rate similar to $T_n(1)$ when $\theta\approx 20^\circ$ on $\mathbb{S}^2$ and $\theta\approx 45^\circ$ on $\mathbb{S}^3$. Regarding the test $T_n^{(10)}$, the simulations show that the parameter $a$ is adequately selected so that it achieves a high power as long as $T_n(1)$ does, while maintaining Type I error. It is interesting to note that for $\theta=135^\circ$ ($\theta=90^\circ$), all (most) tests reject uniformity at a rate significantly lower than $\alpha$, which is attributable to the hyperuniformity created in the generated data.

\begin{table}[ht]
\centering
\small
\begin{tabular}{lR{1.2cm}R{1.2cm}R{1.2cm}R{1.2cm}R{1.2cm}R{1.2cm}R{1.2cm}R{1.2cm}}
    \toprule
    $\theta({}^{\circ})$ & \multicolumn{1}{c}{\small Rayleigh} & \multicolumn{1}{c}{\small Bingham} & \multicolumn{1}{c}{\small Softmax$^{(10)}$} & \multicolumn{1}{c}{\small Poisson$^{(10)}$} & \multicolumn{1}{c}{\small $T_n^{(10)}$} & \multicolumn{1}{c}{\small $T_n(-1)$} & \multicolumn{1}{c}{\small $T_n(0)$} & \multicolumn{1}{c}{\small $T_n(1)$}\\ 
    \hline
    \multicolumn{9}{c}{$q=2$}\\ 
    $1$ & $0.0$ & $35.4$ & $6.0$ & $7.0$ & $\textbf{100.0}$ & $0.0$ & $5.5$ & $\textbf{100.0}$ \\ 
    $10$ & $0.0$ & $34.6$ & $5.0$ & $4.0$ & $54.7$ & $0.1$ & $3.8$ & $\textbf{79.8}$ \\ 
    $20$ & $0.0$ & $\textbf{32.0}$ & $3.5$ & $3.4$ & $12.0$ & $0.4$ & $3.3$ & $\textbf{30.1}$ \\ 
    $45$ & $0.0$ & $\textbf{22.3}$ & $2.3$ & $2.5$ & $3.0$ & $0.8$ & $3.1$ & $9.4$ \\ 
    $90$ & $0.1$ & $4.9$ & $1.7$ & $2.3$ & $2.4$ & $1.9$ & $3.1$ & $4.8$ \\ 
    $135$ & $2.9$ & $2.9$ & $3.0$ & $3.4$ & $3.1$ & $3.7$ & $4.2$ & $4.6$ \\ 
    $180$ & $4.8$ & $4.9$ & $4.4$ & $4.4$ & $3.7$ & $5.0$ & $5.0$ & $5.3$ \\ 
    \hline
    \multicolumn{9}{c}{$q=3$}\\ 
    $1$ & $0.0$ & $49.1$ & $7.9$ & $5.9$ & $\textbf{100.0}$ & $0.0$ & $2.0$ & $\textbf{100.0}$ \\ 
    $10$ & $0.0$ & $48.6$ & $6.5$ & $4.8$ & $\textbf{100.0}$ & $0.0$ & $1.6$ & $\textbf{100.0}$ \\ 
    $20$ & $0.0$ & $45.8$ & $4.7$ & $3.9$ & $75.0$ & $0.0$ & $1.4$ & $\textbf{86.9}$ \\ 
    $45$ & $0.0$ & $\textbf{30.5}$ & $3.0$ & $3.2$ & $13.7$ & $0.0$ & $0.8$ & $23.8$ \\ 
    $90$ & $0.3$ & $5.1$ & $2.2$ & $2.7$ & $2.9$ & $0.6$ & $1.3$ & $5.5$ \\ 
    $135$ & $3.5$ & $3.6$ & $3.4$ & $3.6$ & $3.2$ & $3.4$ & $3.2$ & $3.9$ \\ 
    $180$ & $5.0$ & $4.8$ & $4.8$ & $4.6$ & $4.5$ & $5.1$ & $5.1$ & $5.0$ \\ 
    \bottomrule
\end{tabular}
\caption{\small\label{tbl:unif-antipodal}Rejection percentages for the explored tests on a UAD sample on $\Sq$ of size $n=100$ of angle $\theta$. Note that $\theta=180^{\circ}$ generates an iid sample uniformly distributed on $\Sq$. Boldface indicates the tests whose empirical powers are not significantly smaller than the largest empirical power above the significance level for each row. The significance level is $5\%$.}
\end{table}

\section*{Acknowledgments}

The authors acknowledge support from grant PID2021-124051NB-I00, funded by MCIN/\-AEI/\-10.13039/\-501100011033 and by ``ERDF A way of making Europe''. The second author also acknowledges support from ``Convocatoria de la Universidad Carlos III de Madrid de Ayudas para la recualificación del sistema universitario español para 2021--2023'', funded by Spain's Ministerio de Ciencia, Innovación y Universidades. Comments by two anonymous referees and an Associate Editor are greatly acknowledged.

\appendix

\section{Proofs}\label{sec:proofs}

\begin{proof}[Proof of Proposition \ref{prp:sob_coef}.]
Let $\psi_{\cot}(\theta):=\cot(\theta/2)=(1-\cos\theta)^{-1/2}(1+\cos\theta)^{1/2}$ and $\psi_{\tan}(\theta):=\tan(\theta/2)=(1-\cos\theta)^{1/2}(1+\cos\theta)^{-1/2}$. From Equation ET II 281(4) in \cite{Zwillinger2014} or Equation 16.3.(4) in \citet[p. 281]{Erdelyi1954}, and plugging \eqref{eq:c_d_kq} into \eqref{eq:sob_coef}, we have
\begin{align}
b_{k,q}(\psi_{\cot}) =&\; \frac{1}{c_{k,q}}\,\int_{-1}^1 (1-x)^{-1/2}(1+x)^{1/2} C_{k}^{(q-1)/2}(x) (1-x^2)^{q/2-1}  \, \mathrm{d} x\nonumber\\
=&\; \frac{(2k+q-1)\Gamma((q-1)/2)^2}{2\Gamma(q/2)^2} \, { }_{3} F_{2}\left({\genfrac{}{}{0pt}{}{-k, k+q-1, (q-1)/2}{q/2, q}} ; 1\right)\label{eq:sob_cot},
\end{align}
and, similarly,
\begin{align}
b_{k,q}(\psi_{\tan})=&\; \frac{(2k+q-1)\Gamma((q-1)/2)^2}{2\Gamma(q/2)^2} \, { }_{3} F_{2}\left({\genfrac{}{}{0pt}{}{-k, k+q-1, (q+1)/2}{q/2, q}} ; 1\right)\label{eq:sob_tan}.
\end{align}

Setting $k$ to $2k'$ and $2k'+1$ to distinguish even and odd coefficients, we aim to prove the following equalities:
\begin{align}
{ }_{3} F_{2}\left({\genfrac{}{}{0pt}{}{-2k', 2k'+q-1, (q\mp 1)/2}{q/2, q}}; 1\right)
&=\,\frac{\Gamma \left(k'+\frac{1}{2}\right)^2\Gamma(q/2)^2}{\pi\Gamma \left(k'+q/2\right)^2},\label{eq:3F2_cot_even}\\
{ }_{3} F_{2}\left({\genfrac{}{}{0pt}{}{-2k'-1, 2k'+q, (q\mp 1)/2}{q/2, q}}; 1\right)
&=\pm\frac{2k'+1}{2k'+q}\,\frac{\Gamma \left(k'+1/2\right)^2\Gamma(q/2)^2}{\pi\Gamma \left(k'+q/2\right)^2}.\label{eq:3F2_cot_odd}
\end{align}
Corollary 2.7 in \cite{Lewanowicz1997} gives a generalized Watson's summation formula for\break ${ }_{3} F_{2}\left(a,b,c; \frac{1}{2}(a+b+i+1), 2c+j; 1\right)$ for $i,j\in\mathbb{Z}$ through the following expression, where $a=-2k'-p$, $b=q+2k'+p-1$, and $p=0,1$:
\begin{align*}
    { }_{3} F_{2}&\left({\genfrac{}{}{0pt}{}{-2k'-p, q+2k'+p-1, c}{\frac{1}{2}(q+i), 2c+j}}; 1\right)\\
    &=2^{|i|}\left((q-1)/2\right)_{\mu(1-p)}\left((q+p-1)/2+k'\right)_{\mu p}\\
    &\;\;\;\;\times\frac{\left((q+\mu)/2\right)_{\lfloor|i|/2\rfloor}\left(1/2\right)_{k'+p}\left((q+\mu)/2-c\right)_{k'}}{\left((q+\mu)/2-c\right)_{\lfloor|i|/2\rfloor}\left(c+1/2\right)_{k'+p}\left((q+p-|\mu-p|)/2\right)_{k'}}T_{i,j}^{(p)},
\end{align*}
where $\left(x\right)_{k}:=\Gamma\left(x+k\right)/\Gamma\left(x\right)$ is the Pochhammer symbol, $\mu:=|i|\mod{2}$, and $T_{i,j}^{(p)}=A_{i,j}Q_i^{(p)}+B_{i,j}R_{i}^{(1-p)}$
follows from Equations (2.7), (2.14), and (2.17) in \cite{Lewanowicz1997}.

Equality \eqref{eq:3F2_cot_even} for $c=(q-1)/2$ is proved by plugging in the values $i=p=0$ and $j=1$, where $T_{0,1}^{(0)}=Q_0^{(0)}+R_{0}^{(1)}=1$, with $Q_{0}^{(0)}=1$ and $R_{0}^{(1)}=0$ after some calculations.

Analogously: \eqref{eq:3F2_cot_even} for $c=(q+1)/2$ is proved by setting $i=p=0$ and $j=-1$, giving $\smash{T_{0,-1}^{(0)}}=(2k+q)(1-2k)/q$; \eqref{eq:3F2_cot_odd} for $c=(q-1)/2$ follows with $i=0$ and $j=p=1$, resulting in $\smash{T_{0,1}^{(1)}}=1$; and \eqref{eq:3F2_cot_odd} for $c=(q+1)/2$ results with $i=0$, $j=-1$, and $p=1$, arriving at $\smash{T_{0,-1}^{(1)}}=(2k-1)(2k+q+2)/q$.

Expressions \eqref{eq:sob_weight_odd}--\eqref{eq:sob_weight_even} follow from \eqref{eq:sob_cot}--\eqref{eq:3F2_cot_odd} and the form of $\psi(\cdot;a)$.
\end{proof}

\begin{proof}[Proof of Theorem \ref{thm:asymp_distr}]
    For the Sobolev statistic $S_n(\{w_{k,q}\})$, it follows from Theorem 4.1 and Equation (4.2) in \cite{Gine1975} that
    \begin{align}
        S_n(\{w_{k,q}\}) \inlaw \sum_{k=1}^\infty w_{k,q} Y_k \label{eq:asympSnq}
    \end{align}
    under $\Hcal_0$ as $n\to\infty$. Recall that $\phi(\theta)=\sum_{k=1}^\infty (1+2k/(q-1)) C_k^{(q-1)/2}(\cos\theta)$ almost everywhere for $\theta\in[0,\pi]$ and, hence, $\phi(0)
    =\sum_{k=1}^\infty w_{k,q} d_{k,q}$ from \eqref{eq:c_d_kq}. That is, 
    $S_n(\{w_{k,q}\})=\sum_{k=1}^\infty w_{k,q} d_{k,q}+\frac{2}{n}\sum_{1\leq i<j\leq n} \phi(\theta_{ij})$. The condition $\phi(0)=\sum_{k=1}^\infty w_{k,q} d_{k,q}<\infty$, required for \eqref{eq:asympSnq}, allows us to rewrite \eqref{eq:asympSnq} as
    \begin{align}
        S_n(\{w_{k,q}\}) \inlaw \phi(0)+\sum_{k=1}^\infty w_{k,q} (Y_k-d_{k,q}), \label{eq:Snpsi0}
    \end{align}
    which connects to the result in Theorem 4.3.2 in \cite{Koroljuk1994} after considering the spherical harmonics on $L^2(\mathbb{S}^q,\nu)$ as the orthonormal system for projecting the kernel $(\bx,\by)\in(\mathbb{S}^q)^2\mapsto \phi(\cos^{-1}(\bx'\by))\in\mathbb{R}$.

    The statistic $T_n(a)$ in \eqref{eq:Tna} does not have the ``canonical'' Sobolev $V$-statistic structure because $\psi(0;a)=\infty$ for all $a\in[-1,1]$, and hence the bias term of the diagonal is excluded in its definition. It can be rewritten as
    \begin{align}
        T_{n}(a)=\frac{2}{n}\sum_{1\leq i<j\leq n} \Psi(\bX_i,\bX_j;a)\label{eq:TnaPsi}
    \end{align}
    with $\Psi(\bx,\by;a):=\psi(\cos^{-1}(\bx'\by);a)-\mathbb{E}_{\Hcal_0}\left[\psi(\theta_{12};a)\right]$. This kernel is such that $\mathbb{E}_{\Hcal_0}[\Psi(\bX_1,\bx)]
    =0$ for any $\bx\in\mathbb{S}^q$ (degeneracy) and $\mathbb{E}_{\Hcal_0}[\Psi^2(\bX_1,\bX_2)]<\infty$. The latter follows from $\mathbb{E}_{\Hcal_0}[\Psi^2(\bX_1,\bX_2)]\leq 2(\|\psi(\cdot;a)\|^2_{L_q^2}+b_{0,q}^2(\psi(\cdot;a)))$ and
    \begin{align*}
        \|\psi(\cdot;a)\|^2_{L_q^2}
        &\leq 2\left(\int_{-1}^{1} (1-x)^{q/2-2}\,(1+x)^{q/2}\,\mathrm{d}x + a^2\int_{-1}^{1} (1-x)^{q/2}\,(1+x)^{q/2-2}\,\mathrm{d}x\right)\\
        &\leq 2^{q/2 + 1}\left(\int_{-1}^{1} (1-x)^{q/2-2}\,\mathrm{d}x + a^2\int_{-1}^{1} (1+x)^{q/2-2}\,\mathrm{d}x\right).
    \end{align*}
    If $q>2$, both integrals are finite.

    Therefore, Theorem 4.3.1 in \cite{Koroljuk1994} is readily applicable to \eqref{eq:TnaPsi} and \eqref{eq:null_asymp_dist} follows.
\end{proof}

\begin{proof}[Proof of Theorem \ref{thm:power}]
    Under $\mathrm{P}_{\kappa_n,f}^{(n)}$ and as $n\to\infty$, Theorem 2 in \cite{Garcia-Portugues:Sobolev} gives that
    \begin{align*}
        S_{n,K}(\{w_{k,q}\})\inlaw w_{k_v,q}\widetilde{Y}_{k_v}+\sum_{\substack{k=1\\k\neq k_v}}^{K}w_{k,q}\, Y_k
    \end{align*}
    for a finite Sobolev $V$-statistic $S_{n,K}(\{w_{k,q}\})$ with truncated kernel 
    $\phi_K(\theta)=
    \sum_{k=1}^K (1+2k/(q-1)) w_{k,q}C_k^{(q-1)/2}(\cos\theta)$. In view of the rearrangement done in \eqref{eq:Snpsi0}, this result is readily applicable to $T_{n,K}(a)$ (and any other Sobolev $U$-statistic) given that $\psi_K(0;a)<\infty$, which results in \eqref{eq:nonnull_asymp_dist}.
\end{proof}


\fi

\end{document}